\documentclass{ifacconf}
\usepackage{graphics,epstopdf} % for pdf, bitmapped graphics files
\usepackage{float}
\usepackage{setspace}
\usepackage{mathtools}
\usepackage{graphicx}      % include this line if your document contains figures
\usepackage{natbib}        % required for bibliography
\usepackage{overpic,relsize, here}
\usepackage{empheq}
\usepackage{amssymb}
\theoremstyle{plain}% Theorem-like structures
\newtheorem{theorem}{Theorem}

\newtheorem{lemma}{Lemma}

\theoremstyle{definition}
\newtheorem{definition}{Definition}

\theoremstyle{remark}
\newtheorem{remark}{Remark}

\begin{document}
\begin{frontmatter}

\title{Finite-dimensional output feedback regulator for a mono-tubular heat exchanger process }

%\thanks[footnoteinfo]{}

\author[]{Xiaodong Xu},
\author[UofA]{Stevan Dubljevic}
\address[UofA]{Department of Chemical $\&$ Materials Engineering, University
of Alberta, Canada, T6G 2V4, Stevan.Dubljevic@ualberta.ca}

\begin{abstract}                % Abstract of not more than 250 words.
In this work, we consider the output tracking and disturbance rejection problems of a mono-tubular heat exchanger process and a novel finite-dimensional output feedback regulator is developed. In the proposed output feedback regulator design, measurements available for the regulator do not belong to the set of controlled outputs. In other words, design emphasizes that other than controlled output is used as input signal to the regulator. The proposed output feedback regulator with only plant measurement $y_m(t)$ can realize the exosystem state estimation, disturbance rejection and reference signal tracking, simultaneously. Finally, the results are demonstrated by a mono-tubular heat exchanger process with specific parameters via the computer simulation.
\end{abstract}


\end{frontmatter}
%===============================================================================

\section{Introduction}
The dynamics of many relevant chemical, petrochemical, pharmaceutical  processes can be described by transport-reactions models that take the mathematical form given by partial differential equations (PDEs)  (see \cite{ray1981advanced}, \cite{laabissi2001trajectory}, \cite{christofides2001nonlinear}, \cite{winkin2000dynamical}). The majority of control related research efforts dealing with PDEs have been based on control methods for distributed parameter systems (DPS) which inevitable include the infinite-dimensional nature of these systems (see \cite{curtain1995introduction}, \cite{bensoussan2007representation}, \cite{alizadeh2013boundary}, \cite{deutscher2013finite}). The complexity of infinite dimensional system setting is limiting factor for various controller design and realizations, and in particular needs to be addressed in the framework of fundamental and essential control problems such as problem of stabilization and/or servo-problem.

One classical control problem is the output regulation problem or servo-problem. The problem is formulated as design of controller for fixed plant such that the output can track desired reference signal (and/or reject disturbance) generated by an exosystem. In the past, significant efforts have been made to solve regulator design problems for infinite-dimensional systems: the seminal geometric methods developed in \cite{francis1977linear} in finite-dimensional systems were extended to infinite-dimensional systems (see \cite{byrnes2000output}, \cite{natarajan2014state}). Recently, the generalized geometric methods were introduced for the first order hyperbolic systems in \cite{xu2015output}. Along the same line, the robust regulation problem for infinite-dimensional systems driven by infinite-dimensional exosystems was studied by Pohjolainen in \cite{pohjolainen1982robust} and more recently in \cite{hamalainen2010robust}, \cite{paunonen2010internal}, \cite{paunonen2014robust}. Moreover, for the class of non-spectral systems, Ole Morten Aamo generalised the results based on the backstepping method to deal with boundary disturbances rejection problems and arbitrary-point disturbance rejection problems for linear $2\times 2$ hyperbolic systems (see \cite{aamo2013disturbance}; \cite{aamo2015disturbance}).

In general, there are mainly two versions of the regulator problems: state feedback regulator problem and error feedback regulator problem. In the first problem, the full state information of the plant and exosystem is available to the regulator, while in the second problem, the error signal is measured for the regulator (see e.g., \cite{byrnes2000output} and the references therein), i.e., the regulator design utilizes the knowledge of the plant controlled output and reference signal. Recently, \cite{Deutscher2011output} introduced the finite-dimensional output feedback regulator solving the regulation problem for \emph{Riesz-spectral} system. In the regulator design given by \cite{Deutscher2011output} the measurement output and the controlled output are non-collocated, but it is necessary that the plant measurement and the reference signal genrated by the exosystem are given and available for the regulator design.

In this work, we consider the output regulation problem of a mono-tubular heat exchanger process. A finite-dimensional output feedback regulator design, where only the plant measurement $y_m(t)$ is available for the regulator design and the measurement $y_m(t)$ is different from the controlled output $y(t)$. In particular, the regulator parameters can be easily configured by applying the separation principle. This work is organized to provide the sufficient conditions for the solvability of the output regulation problem. In Section 2, the mono-tubular heat exchanger process, the exosystem and the control objective as well are introduced. In Section 3, based on the feedforward regulator design, the proposed output feedback regulator design utilizes the solution of an auxiliary Sylvester equation. In Section 4, the process parameters are specified to verify the main results. Finally, the remark is made in Section 5.

\section{Problem formulation}
We shall consider the following type of mono-tubular heat exchanger equation on domain $\left\{ {t \in {\mathbb{R}_+ },z \in [0,1]} \right\}$:
\begin{subequations}\label{plant}
\begin{equation}\label{plant-1}
\frac{{\partial x(z,t)}}{{\partial t}} = -\frac{{\partial x(z,t)}}{{\partial z}} + g(z)x(z,t) + b(z)u(t)\vspace{0mm}
\end{equation}
\begin{equation}\label{plant-2}
x(0,t) = 0,x(z,0) = {x_0} \in H\vspace{0mm}
\end{equation}
\begin{equation}\label{plant-3}
y(t) = x({z_1},t),\hspace{1mm}{y_m}(t) = x({z_0},t)
\end{equation}
\end{subequations}
where $x(z,t)\in \mathbb{C}$ is the temperature variation at the time $t$ and at the point $x\in [0,1]$ with respect to an equilibrium point. $u(t)\in \mathbb{C}$ is the control input. $y(t), y_m(t)\in \mathbb{C}$ are the controlled output and the measurement, respectively. $g(z):[0,1] \mapsto \mathbb{R}:z \mapsto g(z)$ is a bounded measurable function, i.e. $g(z)\in L^{\infty}(0,1)$ and $b(z):=\gamma e^{-\bar bz}$ denotes the actuator influence function, where $\bar b, \gamma$ are strictly positive constants.

The equivalent linear infinite-dimensional representation of the system (\ref{plant}) on the state space $H=L^2(0,1)$ is given by \vspace{-5mm}
\begin{eqnarray}
 \dot x(t) &=& Ax(t) + Bu(t),x(0)\in H\label{abplant-1}\\
 y(t) &=& Cx(t), t\ge0\label{abplant-2}\\
 y_m(t) &=& C_mx(t), t\ge0\label{abplant-3}
\end{eqnarray}
with the linear system operator $A$ defined on the domain:
\[D(A) = \left\{ {x \in {L^2}(0,1):{\rm{ }}x(0) = 0} \right\}\]
as:
\begin{equation}\label{Aop}
A = -\frac{d}{{dz}} + g(z) \cdot I
\end{equation}
The input operator is given by $B=b(z)\cdot I$, where $I$ is the identity operator and the controlled output operator is given by $Cx(\cdot,t)=x(z_1,t)$ and the measured output operator is expressed by $C_mx(\cdot,t)=x(z_0,t)$. Here $u(t)$ is element of the input space $U\subset \mathbb{C}$. $y(t),y_m(t)$ are elements of the output spaces $Y,Y_m\subset \mathbb{C}$.

\begin{definition}\label{definition-1}
A linear operator $C$ is called $A-$bounded if $D(A)\subset D(C)$ and if there exist nonnegative constants $a$ and $b$ such that, for all $x(t)\in D(A)$,
\[\left\| {Cx} \right\| \le a\left\| {Ax} \right\| + b\left\| x \right\|\]
\end{definition}
\begin{definition}\label{definition-2}
A linear operator $C$ is called $A-$admissible if it is $A-$bounded and for some $t_f>0$ and $x\in D(A)$ there is a positive constant $M_{t_f}$ such that
\[\int_0^{{t_f}} {{{\left| {C{T_A}(t)x} \right|}^2}} dt \le M_{{t_f}}^2{\left\| x \right\|^2}\]
\end{definition}

From \cite{aksikas2009lq}, we know that $A$ is exponentially stable. From the definition of the input and output operators $B$, $C$ and $C_m$ above, we know that $B\in {\cal L}(U,H)$, $C\in {\cal L}(H_1,Y)$ and $C_m\in {\cal L}(H_1,Y_m)$, i.e. $B$ is bounded and $C$, $C_m$ are unbounded on $H$. The space $H_1$ is a smaller space: $D(A)$ with the norm $\|x\|_1=\|(\beta I-A)x\|$, where $\beta\in \rho(A)$, i.e. resolvent set of $A$. According to \cite{tucsnak2014well}, the operator $C$ is $A-$admissible and the system (\ref{abplant-1}) and (\ref{abplant-2}) is regular. In the same way, it is easy to prove that $C_m$ is $A-$admissible and the system (\ref{abplant-1}) and (\ref{abplant-3}) is regular. Therefore, the system (\ref{abplant-1})-(\ref{abplant-3}) can be redefined as:
\begin{equation}\label{explant}
\begin{array}{l}
\dot x(t) = Ax(t) + Bu(t),x(0) \in H\\
y(t) = {C_\Lambda }x(t),t \ge 0\\
\hspace{-2.5mm}{y_m}(t) = {C_{m\Lambda }}x(t),t \ge 0
\end{array}
\end{equation}
where $C_{\Lambda}\in {\cal L}(H,Y)$ and $C_{m\Lambda}\in{\cal L}(H,Y_m)$ are extension of the operators $C$ and $C_m$ on $H$, respectively (see \cite{tucsnak2009observation} or Def.5.1 in \cite{tucsnak2014well}), defined by
\[\begin{array}{l}
{C_\Lambda }x = \mathop {\lim }\limits_{\lambda  \to  + \infty } C\lambda {(\lambda I - A)^{ - 1}}x,\hspace{1mm}x \in H\\
\hspace{-2.5mm}{C_{m\Lambda }}x = \mathop {\lim }\limits_{\lambda  \to  + \infty } {C_m}\lambda {(\lambda I - A)^{ - 1}}x,\hspace{1mm}x \in H
\end{array}\]
where $\lambda\in \rho(A)$ and $\rho(A)$ is the resolvent set of $A$.

We now consider the linear plant (\ref{explant}) with disturbance $d(t)$:
\begin{equation}\label{explant-d}
\begin{array}{l}
\dot x(t) = Ax(t) + Bu(t)+B_dd(t),x(0) \in H\\
y(t) = {C_\Lambda }x(t),t \ge 0\\
\hspace{-2.5mm}{y_m}(t) = {C_{m\Lambda }}x(t),t \ge 0\vspace{0mm}
\end{array}
\end{equation}
where the disturbance is $d(t)\in U_d$ and $U_d$ is a real Hilbert space. The disturbance location operator $B_d\in {\cal L}(U_d,H)$ are bounded on $H$. Since the operators $C\in {\cal L}(H_1,Y)$ and $C_m\in {\cal L}(H_1,Y_m)$ are $A-$admissible for $e^{At}$, the system (\ref{explant-d}) is well-posed and for some (hence, for every) $s\in \rho(A)$, ${C_\Lambda }{(sI - A)^{ - 1}}B$, ${C_\Lambda }{(sI - A)^{ - 1}}{B_d}$, ${C_{m\Lambda }}{(sI - A)^{ - 1}}B$ and ${C_{m\Lambda }}{(sI - A)^{ - 1}}{B_d}$ exist.

The disturbance $d(t)$ and the reference trajectory $y_r(t)$ are generated by a known autonomous finite-dimensional signal process (exogenous system) described by\vspace{-1mm}
\begin{equation}\label{exosystem-1}\vspace{0mm}
\hspace{17mm}\dot w(t)= Sw(t), \hspace{1mm}w(0)\in W
\end{equation}
\begin{equation}\label{exosystem-2}\vspace{0mm}
d(t) = Fw(t)
\end{equation}
\begin{equation}\label{exosystem-3}\vspace{0mm}
\hspace{-1mm}y_r(t) = Qw(t)
\end{equation}
where $S$ is a diagonalizable and skew-hermitian matrix having all its eigenvalues on the imaginary axis, i.e. $iw_k$ with $i=\sqrt{-1}$, $F$ and $Q$ are matrices of appropriate dimensions. $w(t)\in W$ is the state of exosystem and the space $W$ is a $n-$dimensional complex Euclidean space, i.e., $W\subset \mathbb{C}^n$. Therefore, we have
\begin{equation}\label{S-seri}
Sw = \sum\limits_{k = 1}^n {i{w_k}\left\langle {w,{\psi _k}} \right\rangle {\psi _k}}
\end{equation}
where ${\left( {{\psi _k}} \right)_{k \in \mathbb{N}}}$ is an orthonormal basis of $\mathbb{C}^n$. By assigning the different values of the eigenvalues of $S$ and appropriate $Q$, the exosystem can generate steplike, harmonic and ramp trajectories.

Let $X_c=H \oplus W$ be the composite state-space, consisting of the states of plant (\ref{explant-d}) and exosystem (\ref{exosystem-1}), and we obtain the composite system\vspace{0mm}
\begin{equation}\label{composite-1}
\dot {x}_c(t) = A_c x_c(t) + B_c u(t),\hspace{2mm}t>0,\hspace{2mm}x_c(0)\in X_c
\end{equation}
\begin{equation}\label{composite-2}\vspace{0mm}
\hspace{-21.5mm}e(t) = C_c x_c(t)
\end{equation}
\begin{equation}\label{composite-2}\vspace{0mm}
\hspace{-24mm}y_m(t) = C_y x_c(t)
\end{equation}
with the composite state $x_c(t)={x_c}(t) = {\left[ {\begin{array}{*{20}{c}}
   {x(t)^T} & {w{{(t)}^T}}  \\
\end{array}} \right]^T}$, where the composite system operator is given by
\[{A_c} = \left[ {\begin{array}{*{20}{c}}
   A & P  \\
   0 & S  \\
\end{array}} \right] \text{ with }P = {B_d}F\]
where $D({A_c}) = D(A) \oplus W \subset {X_c}$, the operators $C_c\in {\cal L}(X_c, Y)$, $C_y\in {\cal L}(X_c, Y_m)$ and $B_c\in {\cal L}(U, X_c)$ are given by
${C_c} = \left[ {\begin{array}{*{20}{c}}
   C_{\Lambda} & { - Q}  \\
\end{array}} \right],\hspace{2mm}{C_y} = \left[ {\begin{array}{*{20}{c}}
   C_{\Lambda m} & { 0}  \\
\end{array}} \right],\hspace{2mm}{B_c} = \left[ {\begin{array}{*{20}{c}}
   B  \\
   0  \\
\end{array}} \right]$.

It is evident that the system generator $A_c$ generates a $C_0-$semigroup on $X_c$.

In this work, we propose to design finite-dimensional output feedback regulator such that
\begin{equation}\label{error}
\mathop {\lim }\limits_{t \to \infty } \left( {y(t) - {y_r}(t)} \right) = 0,\hspace{1mm}\forall x(0) \in H,\hspace{1mm}w(0) \in W
\end{equation}
in the sense of exponential tracking error $(e(t)=y(t)-y_r(t))$ decay.

\section{Output feedback regulator problem}

We now first demonstrate the design of the finite-dimensional feedforward regulator on the space $W$:
\begin{equation}\label{sim-for-regulator}
\begin{array}{l}
 {{\dot r}_w}(t) = S{r_w}(t) \\
\hspace{2mm} u(t) = {\Gamma}{r_w}(t) \\
 \end{array}
\end{equation}
with initial condition $r_w(0)=r_{w0}\in W$ and $\Gamma\in {\cal L}(W,U)$. By applying the feedforward regulator, the composite system (\ref{composite-1}) tracks the linear function
\begin{equation}\label{lin-fun}
x_c(t) = \tilde \Pi r_w(t)
\end{equation}
where $\tilde \Pi  = \left[ {\begin{array}{*{20}{c}}
   \Pi   \\
   I  \\
\end{array}} \right]$ with the bounded operator $\Pi\in {\cal L}(W,H)$ with $\Pi D(S)\subset D(A)$, if the initial condition is  given by
\begin{equation}\label{in-codi}
{x_c}(0) = \tilde \Pi {r_w}(0)
\end{equation}

The following theorem provides conditions such that (\ref{sim-for-regulator}) solves output regulation problem.

\begin{theorem}\label{theorem-2}
For the system (\ref{explant-d}), the feedforward regulator (\ref{sim-for-regulator}) achieves exponentially stable tracking in the sense of (\ref{error}) for initial values (\ref{in-codi}), if there exist operators $\Pi\in {\cal L}(W,H)$ defined in (\ref{lin-fun}) and $\Gamma\in {\cal L}(W,U)$ satisfying the Sylvester equations:\vspace{-2mm}
\begin{equation}\label{Syl-eq-1}
\hspace{11mm} \Pi S - A\Pi  = B\Gamma  + P \vspace{0mm}
 \end{equation}
 \begin{equation}\label{Syl-eq-2}\vspace{0mm}
 C_{\Lambda}\Pi  - Q   = 0
\end{equation}
where the eigenvalues of $S$ do not coincide with an invariant zero of the plant (\ref{abplant-1})-(\ref{abplant-3}).\\
In this case, the feedforward regulator (\ref{sim-for-regulator}) has the form:\vspace{-1mm}
\begin{equation}\label{reg-new-form}
\begin{array}{l}
 {{\dot r}_w}(t) = S{r_w}(t) \\
\hspace{2mm} u(t) = \Gamma{r_w}(t) \vspace{-3mm}\\
 \end{array}
\end{equation}
\end{theorem}
\textbf{Proof.} The dynamics of $x_c(t)-\tilde \Pi r_w(t)$ are described by the autonomous abstract differential equation ${{\dot x}_c}(t)- \tilde \Pi \dot r_w(t) = {A_c}[{x_c}(t) - \tilde \Pi {r_w}(t)]$ with the initial value ${x_c}(0) - \tilde \Pi {r_w}(0)\in X_c$ if there exists a bounded linear operator $\tilde \Pi$ such that the following Sylvester operator equation holds:\vspace{-1mm}
 \begin{equation}\label{tem-syl}
\left[ {\begin{array}{*{20}{c}}
   A & P  \\
   0 & S  \\
\end{array}} \right]\tilde \Pi  - \tilde \Pi S =- \left[ {\begin{array}{*{20}{c}}
   B  \\
   0  \\
\end{array}} \right]\Gamma
 \end{equation}

Since $A_c$ is the generator of an infinitesimal $C_0-$semigroup, the initial value problem has a unique solution. If (\ref{in-codi}) holds, the solution of ${{\dot x}_c}(t)- \tilde \Pi \dot r_w(t) = {A_c}[{x_c}(t) - \tilde \Pi {r_w}(t)]$ is ${x_c}(t) - \tilde \Pi {r_w}(t)=0$. Then, the plant state $x(t)$ and the state of exosystem $w(t)$ can be expressed by $x(t)=\Pi r_w(t)$ and $w(t)=r_w(t)$, so that the tracking error in (\ref{error}) takes the form $y(t)-y_r(t)=C_{\Lambda}\Pi r_w(t)-Qr_w(t)$. Therefore, in order to achieve output regulation (\ref{error}), the operator $\Pi$ has in addition to satisfy $C_{\Lambda}\Pi-Q=0$ since $S$ is anti-Hurwitz matrix. A direct calculation shows that the equation (\ref{Syl-eq-1}) is equivalent to the equation (\ref{tem-syl}). In \cite{byrnes2000output}, it is shown that there exists a solution $(\Pi, \Gamma)$ of the Sylvester equation (\ref{Syl-eq-1})-(\ref{Syl-eq-2}), if the eigenvalues of $S$ do note coincide with an invariant zero of the plant (\ref{abplant-1})-(\ref{abplant-3}).\qed

Since in general case the state of exosystem cannot be measured, the initial values $r_w(0)$ of feedforward controller (\ref{sim-for-regulator}) cannot be chosen such that $(\ref{in-codi})$ holds. Therefore, in general the initial error $x_c(0)-\Pi r_w(0)\ne 0$ may result so that output regulation is not achieved. However, the exponentially decaying tracking error in (\ref{error}) can be obtained by stabilizing the dynamics of $x_c(t)-\tilde \Pi r_w(t)$. Therefore, in this section, we propose a output feedback regualtor. Following internal model principle, we design the following regulator by including the internal model of exosystem:\vspace{-3mm}
\begin{empheq}{align}
{{\dot r}_w}(t) &= S{r_w}(t) + {L_y}\left[ {{y_m}(t) - {C_y}\tilde \Pi {r_w}(t)} \right] \label{ym-contr-1}\\
u(t) &= \Gamma {r_w}(t) + k\left[ {{y_m}(t) - {C_y}\tilde \Pi {r_w}(t)} \right]\label{ym-contr-2}
\end{empheq}
with ${C_y} = \left[ {\begin{array}{*{20}{c}}
{{C_{\Lambda m}}}&0
\end{array}} \right]$.\\
The output feedback regulator (\ref{ym-contr-1})-(\ref{ym-contr-2}) is an extension of the feedforward regulator (\ref{sim-for-regulator}). First of all, we investigate the dynamics of $x_c(t)-\tilde \Pi r_w(t)$ as the following:
\begin{equation}\label{tem-tr}
\begin{array}{l}
{{\dot x}_c}(t) - \tilde \Pi {{\dot r}_w}(t) = \left( {{A_c} + k{B_c}{C_y} - \tilde \Pi {L_y}{C_y}} \right){x_c}(t)\\
 - \left( {\tilde \Pi S - {B_c}\Gamma  + k{B_c}{C_y}\tilde \Pi  - \tilde \Pi {L_y}{C_y}\tilde \Pi } \right){r_w}(t)\\
 = \left( {{A_c} + k{B_c}{C_y} - \tilde \Pi {L_y}{C_y}} \right){x_c}(t)\\
 - \left( {{A_c}\tilde \Pi  + k{B_c}{C_y}\tilde \Pi  - \tilde \Pi {L_y}{C_y}\tilde \Pi } \right){r_w}(t)\\
 = \left( {{A_c} + k{B_c}{C_y} - \tilde \Pi {L_y}{C_y}} \right)\left( {{x_c}(t) - \tilde \Pi {r_w}(t)} \right)\\
 = {{\hat A}_c}\left( {{x_c}(t) - \tilde \Pi {r_w}(t)} \right)
\end{array}
\end{equation}
where we have used the equation (\ref{tem-syl}). From Theorem 1 and (\ref{tem-tr}), it is easy to conclude the following lemma.
\begin{lemma}\label{lemma-1}
For the plant (\ref{abplant-1})-(\ref{abplant-3}) and the exosystem (\ref{exosystem-1})-(\ref{exosystem-3}) defined in Section 2, the finite-dimensional regulator (\ref{ym-contr-1})-(\ref{ym-contr-2}) solves the output regulation problem (\ref{error}) if the operators $\Pi\in {\cal L}(W,H)$ with $\Pi D(S)\subset D(A)$ and $\Gamma\in {\cal L}(W,U)$ satisfy the Sylvester equations (\ref{Syl-eq-1})-(\ref{Syl-eq-2}) and if there exist controller parameters $L_y$ and $k$ such that the operator $\hat A_c$ in (\ref{tem-tr}) is the infinitesimal generator of an exponentially stable $C_0-$semigroup.
\end{lemma}

Consequently, we investigate the choice of the regulator parameters $L_y$ and $k$ such that the operator $\hat A_c$ generates an exponentially stable $C_0-$semigroup and thus the output regulation problem is solved.
\begin{theorem}\label{theorem-3}
The controller (\ref{ym-contr-1})-(\ref{ym-contr-2}) stabilizes the operator $\hat A_c$, provided that $L_y$ and $k$ are chosen as follows: We assume $k=k_1+k_2$. Choose $k_1$ such that the operator $A+k_1BC_{\Lambda m}$ generates an exponentially stable $C_0-$semigroup. $L_y$ can be chosen such that $S+L_yC_{\Lambda m}\Pi_0$ is exponentially stable, where $\Pi_0\in {\cal L}(W,H)$ is the following solution of the following auxiliary Sylvester equation
\begin{equation}\label{sylv-3}
{\Pi _0}S -( A + {k_1}B{C_{\Lambda m}}){\Pi _0} = - P
\end{equation}
with $\Pi_0 D(S)\subset D(A)$. Finally, we choose $k_2=\frac{\Pi L_y+\Pi_0 L_y}{b(z)}$, since $b(z)$ is nonzero spatially varying function.
\end{theorem}
\textbf{Proof:} With the bounded similarity transformations, the operator $\hat A_c$ is transformed into block lower triangular form, where diagonal blocks can be stabilized by choosing appropriate regulator parameters. by assuming $k=k_1+k_2$, we rewrite $\hat A_c$ as:
\[\begin{array}{l}
{{\hat A}_c} = \left[ {\begin{array}{*{20}{c}}
A&P\\
0&S
\end{array}} \right] + \left[ {\begin{array}{*{20}{c}}
{kB - \Pi {L_y}}\\
{ - {L_y}}
\end{array}} \right]\left[ {\begin{array}{*{20}{c}}
{{C_{\Lambda m}}}&0
\end{array}} \right]\vspace{1mm}\\
 = \left[ {\begin{array}{*{20}{c}}
A&P\\
0&S
\end{array}} \right] + \left[ {\begin{array}{*{20}{c}}
{{k_1}B}\\
0
\end{array}} \right]\left[ {\begin{array}{*{20}{c}}
{{C_{\Lambda m}}}&0
\end{array}} \right]
\end{array}\]
\[\begin{array}{l}
 + \left[ {\begin{array}{*{20}{c}}
{{k_2}B}\\
0
\end{array}} \right]\left[ {\begin{array}{*{20}{c}}
{{C_{\Lambda m}}}&0
\end{array}} \right] + \left[ {\begin{array}{*{20}{c}}
{ - \Pi {L_y}}\\
{ - {L_y}}
\end{array}} \right]\left[ {\begin{array}{*{20}{c}}
{{C_{\Lambda m}}}&0
\end{array}} \right]\vspace{1mm}\\
 = \left[ {\begin{array}{*{20}{c}}
{A + {k_1}B{C_{\Lambda m}}}&P\\
0&S
\end{array}} \right] + \left[ {\begin{array}{*{20}{c}}
{{k_2}B - \Pi {L_y}}\\
{ - {L_y}}
\end{array}} \right]\left[ {\begin{array}{*{20}{c}}
{{C_{\Lambda m}}}&0
\end{array}} \right]
\end{array}\]

By applying the following similarity transformation:
\[T = \left[ {\begin{array}{*{20}{c}}
   I & {{\Pi _0}}  \\
   0 & I  \\
\end{array}} \right],\hspace{2mm}{T^{ - 1}} = \left[ {\begin{array}{*{20}{c}}
   I & { - {\Pi _0}}  \\
   0 & I  \\
\end{array}} \right]\]
we transform $\hat A_c$ into the form:
\begin{equation}\label{tem-Ac}
\begin{array}{*{20}{l}}
\begin{array}{l}
T{{\hat A}_c}{T^{ - 1}} = \left[ {\begin{array}{*{20}{c}}
I&{{\Pi _0}}\\
0&I
\end{array}} \right]\left[ {\begin{array}{*{20}{c}}
{A + {k_1}B{C_{\Lambda m}}}&P\\
0&S
\end{array}} \right]\left[ {\begin{array}{*{20}{c}}
I&{ - {\Pi _0}}\\
0&I
\end{array}} \right]\vspace{1mm}\\
\hspace{17mm} + \left[ {\begin{array}{*{20}{c}}
I&{{\Pi _0}}\\
0&I
\end{array}} \right]\left[ {\begin{array}{*{20}{c}}
{{k_2}B - \Pi {L_y}}\\
{ - {L_y}}
\end{array}} \right]\left[ {\begin{array}{*{20}{c}}
{{C_{\Lambda m}}}&0
\end{array}} \right]\left[ {\begin{array}{*{20}{c}}
I&{ - {\Pi _0}}\\
0&I
\end{array}} \right]\vspace{1mm}
\end{array}\\
\begin{array}{l}
\hspace{14mm} = \left[ {\begin{array}{*{20}{c}}
{A + {k_1}B{C_{\Lambda m}}}&{P + {\Pi _0}S}\\
0&S
\end{array}} \right]\left[ {\begin{array}{*{20}{c}}
I&{ - {\Pi _0}}\\
0&I
\end{array}} \right]\\
\hspace{16mm} + \left[ {\begin{array}{*{20}{c}}
{{k_2}B - \Pi {L_y} - {\Pi _0}{L_y}}\\
{ - {L_y}}
\end{array}} \right]\left[ {\begin{array}{*{20}{c}}
{{C_{\Lambda m}}}&{ - {C_{\Lambda m}}{\Pi _0}}
\end{array}} \right]\vspace{1mm}
\end{array}\\
{{\textrm{let }}{k_2}B - \Pi {L_y} - {\Pi _0}{L_y} = 0} \vspace{1mm}\\
\hspace{14mm}  = \left[ {\begin{array}{*{20}{c}}
{A + {k_1}B{C_{\Lambda m}}}&{P + {\Pi _0}S - A{\Pi _0} - {k_1}B{C_{\Lambda m}}{\Pi _0}}\\
{ - {L_y}{C_{\Lambda m}}}&{S + {L_y}{C_{\Lambda m}}{\Pi _0}}
\end{array}} \right]
\end{array}
\end{equation}
For specified $k_1$, if we let $\Pi_0$ satisfy the following auxiliary Sylvester equation:
\[{\Pi _0}S - A{\Pi _0} - {k_1}B{C_{\Lambda m}}{\Pi _0} + P = 0\]
with $\Pi_0D(S)\subset D(A)$, then, $\hat A_c$ can be written as block lower triangular form:
\[T{{\hat A}_c}{T^{ - 1}} = \left[ {\begin{array}{*{20}{c}}
{A + {k_1}B{C_{\Lambda m}}}&0\\
{ - {L_y}{C_{\Lambda m}}}&{S + {L_y}{C_{\Lambda m}}{\Pi _0}}
\end{array}} \right]\]

From above equation, we choose $k_1$ and $L_y$ such that $A+k_1BC_{\Lambda m}$ generates an exponentially stable $C_0-$semigroup, $S+L_yC_{\Lambda m}\Pi_0$ is stable and thus $\hat A_c$ generates an exponentially stable $C_0-$semigroup. Once $k_1$, $L_y$ and $\Pi_0$ are obtained, $k_2$ can be easily calculated from:
\[k_2B-\Pi L_y-\Pi_0 L_y=0\]
Since $B=b(z)\cdot I\ne 0$, $k_2=\frac{\Pi L_y+\Pi_0 L_y}{B}$. Therefore, we can choose $k=k_1+\frac{\Pi L_y+\Pi_0 L_y}{b(z)}$.

\begin{remark}\label{remark-2}
It can be noticed that for the construction of the output feedback regulator, only the measured output $y_m(t)$ is utilized and the reference signal $y_r(t)$ is not necessary. From previous section, we know that the proposed regulator has the form (\ref{ym-contr-1})-(\ref{ym-contr-2}) because our target is to drive the composite state $x_c(t)$ tracks the function $\tilde \Pi r_w(t)$. More precisely, the mission of the output feedback regulator is to drive $w(t)$ tracking $r_w(t)$ and $x(t)$ tracking $\Pi r_w(t)$, which indicates that $x(t)$ finally tracks $\Pi w(t)$. Therefore, the measured output $y_m(t)=C_{\Lambda m}x(t)$ tracks $C_{\Lambda m}\Pi w(t)$. This implies that $y_m(t)$ contains the state information of the exosystem approximately. Moreover, (\ref{ym-contr-1}) is actually an observer of the exosystem. In particular, since the reference signal $y_r(t)$ is not available, it is necessary that disturbance is nonzero, i.e., $d(t)\ne 0$ and $B_d\ne 0$ in (\ref{explant-d}). \qed\vspace{0mm}
\end{remark}

\section{Numerical simulation}
Let us set $g(z)=0.4z$, $\gamma=1$ and $\bar b=0.5$ in system (\ref{plant}). We assume that the controlled output $y(t)$ is the evolution of the state at point $z_1=0.5$ and the measured output $y_m(t)$ is the evolution of the state at boundary point $z_0=1$. Obviously, the controlled output and the measured output are different.

In this section, our objective is to design a finite-dimensional output feedback regulator such that the controlled output $y(t)$ tracks a given reference signal of the form $y_r(t)=\Upsilon sin(\alpha t)$ and rejects the disturbance $d(t)=\Upsilon\cos(\alpha t)$ as well. The disturbance $d(t)$ and the reference signal $y_r(t)$ are generates by the exosystem if the form of (\ref{exosystem-1})-(\ref{exosystem-3}):
\begin{equation}
\left[ {\begin{array}{*{20}{c}}
{{{\dot w}_1}(t)}\\
{{{\dot w}_2}(t)}
\end{array}} \right] = \left[ {\begin{array}{*{20}{c}}
0&\alpha \\
{ - \alpha }&0
\end{array}} \right]\left[ {\begin{array}{*{20}{c}}
{{w_1}(t)}\\
{{w_2}(t)}
\end{array}} \right],\left[ {\begin{array}{*{20}{c}}
{{w_1}(0)}\\
{{w_2}(0)}
\end{array}} \right] = \left[ {\begin{array}{*{20}{c}}
0\\
\Upsilon
\end{array}} \right]
\end{equation}
We can computer the solution of the exosystem system with the parameter values: $\alpha =2$ and $\Upsilon =5$. In terms of notation defined in (\ref{exosystem-2}) and (\ref{exosystem-3}), we set $Q = \left[ {\begin{array}{*{20}{c}}
1&0
\end{array}} \right]$ and $F = \left[ {\begin{array}{*{20}{c}}
0&{1}
\end{array}} \right]$. Then, the disturbance $d(t)=5cos(2t)$ and the reference signal $y_r(t)=5sin(2t)$. From the form of $S$: $S = \left[ {\begin{array}{*{20}{c}}
0&2\\
{ - 2}&0
\end{array}} \right]$, it is easy to calculate its eigenvalues and orthonormal eigenvectors:\vspace{-1mm}
\[i{w_1} = 2i,\hspace{3mm}i{w_2} =  - 2i \vspace{-1mm}\]
\[{\psi _1} = \left[ {\begin{array}{*{20}{c}}
{ - 0.7071i}\\
{0.7071}
\end{array}} \right],\hspace{2mm}{\psi _2} = \left[ {\begin{array}{*{20}{c}}
{0.7071i}\\
{0.7071}
\end{array}} \right]\]

According to the content in Section 2, we know the plant (\ref{plant}) is regular. The input operator $B\in {\cal L}(U,L^2(0,1))$ and disturbance location operator $B_d\in {\cal L}(U_d,L^2(0,1))$, where $U, U_d\subset \mathbb{R}$, are defined by:
\[B =b(z)I =e^{-0.5z}I,\hspace{2mm}{B_d} = 0.2\]
\begin{figure}[H]
  \centering
  \includegraphics[width=3in]{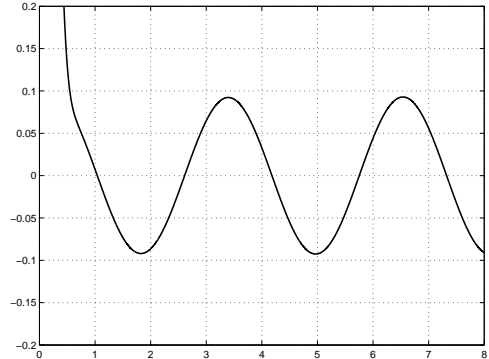}\vspace{-7mm}\\
  \caption{The controlled output $y(t)$ driven by the disturbance $d(t)=5cos(2t)$ with $u(t)=0$.}\label{fig-1-2}
\end{figure}
We now carry on solving the constrained Sylvester equation (\ref{Syl-eq-1})-(\ref{Syl-eq-2}) with $P = {B_d}F = \left[ {\begin{array}{*{20}{c}}
0&{0.2}
\end{array}} \right]$ in Theorem 2. In this example, we construct the operators: $\Gamma  = \left[ {\begin{array}{*{20}{c}}
{{\gamma _1}}&{{\gamma _2}}
\end{array}} \right] \in {\cal L}({\mathbb{C}^2},U)$ and $\Pi  = \left[ {\begin{array}{*{20}{c}}
{{\Pi _1}(z)}&{{\Pi _2}(z)}
\end{array}} \right] \in {\cal L}({\mathbb{C}^2},{L^2}(0,1))$. According to \cite{byrnes2000output}, \cite{natarajan2014state} and \cite{xu2015output}, the solution of constrained Sylvester equation (\ref{Syl-eq-1})-(\ref{Syl-eq-2}) is given by:
\[\Pi {\psi _1} = {\left( {2iI - A} \right)^{ - 1}}(B\Gamma  + P){\psi _1}\]
The straightforward calculation shows that
\[\begin{array}{l}
\hspace{-3mm}\left[ {\begin{array}{*{20}{c}}
{{\Pi _1}(z) - {{\left( {2iI - A} \right)}^{ - 1}}B{\gamma _1},}&{{\Pi _2}(z) - {{\left( {2iI - A} \right)}^{ - 1}}\left( {B{\gamma _2} + 0.2} \right)}
\end{array}} \right]\\
 \times \left[ {\begin{array}{*{20}{c}}
{ - i}\\
1
\end{array}} \right] = 0
\end{array}\]
and we obtain the identity
\begin{equation}\label{Pi-12}
{\Pi _2}(z) - {\Pi _1}(z)i = \left( {{\gamma _2} - {\gamma _1}i} \right){\left( {2iI - A} \right)^{ - 1}}b(z) + 0.2{\left( {2iI - A} \right)^{ - 1}}
\end{equation}
Moreover, the equation (\ref{Syl-eq-2}) reduces to
\[{C_\Lambda }{\Pi _1}(z) = 1,\hspace{3mm}{C_\Lambda }{\Pi _2}(z) = 0\]
Now we apply $C_{\Lambda}$ to the equation (\ref{Pi-12}), then we obtain
\begin{equation}\label{gama-12}
\begin{array}{*{20}{l}}
{ - i = \left( {{\gamma _2} - {\gamma _1}i} \right)G(2i) + 0.2{G_1}(2i)}\\
\begin{array}{l}
 = \left( {{\gamma _2} - {\gamma _1}i} \right)\left( {{\rm{Re}}\left( {G(2i)} \right) + i{\rm{Im}}\left( {G(2i)} \right)} \right)\\
 + 0.2\left( {{\rm{Re}}\left( {{G_1}(2i)} \right) + i{\rm{Im}}\left( {{G_1}(2i)} \right)} \right)
\end{array}
\end{array}
\end{equation}
where $\text{Re}$ and $\text{Im}$ denote the real and imaginary parts and
\[\hspace{-3mm}G(s) = {C_\Lambda }{\left( {sI - A} \right)^{ - 1}}b(z),s \in \rho (A)\]
\[ {G_1}(s) = {C_\Lambda }{\left( {sI - A} \right)^{ - 1}},s \in \rho (A)\]
where $G(s)$ is the system transfer function. We can regard $G_1(s)$ as $G(s)$ with $b(z)=1$.

From (\ref{gama-12}), we can obtain the following equations in terms of $G(2i)$ and $G_1(2i)$:
\[\begin{array}{l}
{\gamma _1}{\mathop{\rm Im}\nolimits} \left( {G(2i)} \right) + {\gamma _2}{\mathop{\rm Re}\nolimits} \left( {G(2i)} \right) =  - 0.2{\mathop{\rm Re}\nolimits} \left( {{G_1}(2i)} \right)\\
{\gamma _1}{\mathop{\rm Re}\nolimits} \left( {G(2i)} \right) - {\gamma _2}{\mathop{\rm Im}\nolimits} \left( {G(2i)} \right) = 1 + 0.2{\mathop{\rm Im}\nolimits} \left( {{G_1}(2i)} \right)
\end{array}\]
These equations allow us to get the explicit expression of $\gamma_1$ and $\gamma_2$:
\[\begin{array}{l}
{\gamma _1} = \frac{{\left( {1 + 0.2{\rm{Im}}\left( {{G_1}(2i)} \right)} \right){\rm{Re}}\left( {G(2i)} \right)}}{{{{\left| {G(2i)} \right|}^2}}}\\
\hspace{9mm} - \frac{{0.2{\rm{Re}}\left( {{G_1}(2i)} \right){\rm{Im}}\left( {G(2i)} \right)}}{{{{\left| {G(2i)} \right|}^2}}}
\end{array}\]
\[\begin{array}{l}
{\gamma _2} =  - \frac{{0.2{\rm{Re}}\left( {{G_1}(2i)} \right){\rm{Re}}\left( {G(2i)} \right)}}{{{{\left| {G(2i)} \right|}^2}}}\\
\hspace{9mm} - \frac{{\left( {1 + 0.2{\rm{Im}}\left( {{G_1}(2i)} \right)} \right){\rm{Im}}\left( {G(2i)} \right)}}{{{{\left| {G(2i)} \right|}^2}}}
\end{array}\]
It is apparent that the specified value of the transfer function is essential to the calculation of $\gamma_1$ and $\gamma_2$.

\begin{remark}\label{remark-3}
In most cases, it is not easy to calculate the explicit expression of transfer function for the PDE systems. However, from the expression of $\gamma_1$ and $\gamma_2$, it is enough to calculate the value of $G(2i)$ and $G_1(2i)$. Therefore, in this work, we propose a numerical method to provide value of $G(2i)$ and $G_1(2i)$. It is well known that if the input $u(t)$ of the system (\ref{plant}) is impulse signal $\delta(t)$, then, the output $y(t)$ is the impulse response whose Laplace transform $Y(s) = {L_t}{\rm{\{ }}y(t)\}(s)$ is actually the transfer function $G(s)$, i.e. $G(s)=Y(s)$. More precisely, in order to calculate $G(2i)$, we just need to calculate $Y(2i)$. \hspace{5mm}$\rule{2mm}{2mm}$
\end{remark}
According to Remark \ref{remark-3}, we substitute impulse input $u(t)=\delta(t)$ into (\ref{plant}) and take Laplace transform on both side of (\ref{plant}) in time by setting zero initial condition,
\begin{equation}\label{lpls}
\begin{array}{l}
\frac{{\partial x(z,s)}}{{\partial z}} = \left( {0.4z - s} \right)x(z,s) + {e^{ - 0.5z}}U(s)\\
x(0,s) = 0,\\
Y(s) = x(0.5,s)
\end{array}
\end{equation}
with $U(s)={\cal L}_t{u(t)}(s)=1$. Our goal is to computer the value of $G(2i)=Y(2i)=x(0.5,2i)$. Therefore, one can directly to solve (\ref{lpls}) with $s=2i$, by numerical method, e.g. finite difference. It is easy to get:
\[G(2i) = x(0.5,2i) = 0.3611 - 0.2021i\]
In the same way, to calculate $G_1(2i)$, one just needs to set $b(z)=1$ in (\ref{lpls}) and obtains,
\[{G_1}(2i) = 0.4121 - 0.2183i\]

With the values of $G(2i)$ and $G_1(2i)$, one can calculate,
\[{\gamma _1} = 2.114,\hspace{3mm}{\gamma _2} = 0.9549\]

We write the equation (\ref{Syl-eq-1}),\vspace{-1mm}
\begin{subequations}\label{pi-tem}
\begin{equation}\label{pi-tem-1}
{{\Pi '}_1}(z) - g(z){\Pi _1}(z) - 2{\Pi _2}(z) = b(z){\gamma _1}\vspace{0mm}
\end{equation}
\begin{equation}\label{pi-tem-2}
{{\Pi '}_2}(z) - g(z){\Pi _2}(z) + 2{\Pi _1}(z) = b(z){\gamma _2} + 0.2\vspace{0mm}
\end{equation}
\begin{equation}\label{pi-tem-3}
{\Pi _1}(0) = 0,{\Pi _2}(0) = 0
\end{equation}
\end{subequations}
The boundary condition (\ref{pi-tem-3}) is necessary such that $\Pi_1$ and $\Pi_2$ lie in the domain of $A$: $D(A)$. Obviously, the equation (\ref{pi-tem}a-b-c) can be solved off-line numerically.

Now we carry on the construction of output feedback regulator of form (\ref{ym-contr-1})-(\ref{ym-contr-2}). Under the control of the regulator (\ref{ym-contr-1})-(\ref{ym-contr-2}), we show that the controlled output $y(t)$ tracks the reference signal $y_r(t)=5sin(2t)$ in spite of the existence of disturbance $d(t)=5cos(2t)$.
We assume that only the measured output of the plant: $y_m(t)$ is available to the regulator and the controlled output $y(t)$ and the measured output $y_m(t)$ are different, i.e., $y(t)=x(0.5,t)$ and $y_m(t)=x(1,t)$.

In this part, the key is to calculate $k$ and $L_y$. According to Theorem \ref{theorem-3}, one can first set $k_1=0$ since $A$ is exponentially stable. Then, by setting ${\Pi _0}(z) = \left[ {\begin{array}{*{20}{c}}
{{\Pi _{01}}(z)}&{{\Pi _{02}}(z)}
\end{array}} \right]$ one can rewrite (\ref{sylv-3}) as:\vspace{-1mm}
\[\begin{array}{l}
{{\Pi '}_{01}}(z) - g(z){\Pi _{01}}(z) - 2{\Pi _{02}}(z) = 0\\
{{\Pi '}_{02}}(z) - g(z){\Pi _{02}}(z) + 2{\Pi _{01}}(z) =  - 0.2\\
{\Pi _{01}}(0) = 0,{\Pi _{02}}(0) = 0
\end{array}\]

\begin{figure}[ht]
  \centering
  \includegraphics[width=3in]{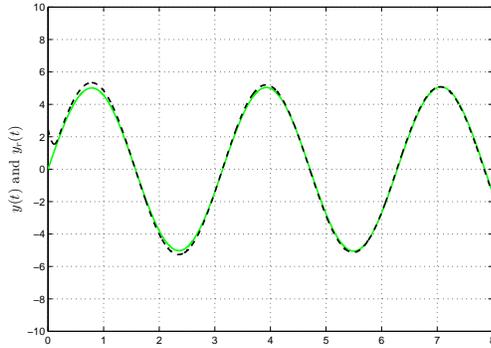}\vspace{-7mm}\\
  \caption{The controlled output $y(t)$ tracks the reference signal $y_r(t)=5sin(2t)$ under the control of regulator (\ref{ym-contr-1})-(\ref{ym-contr-2})}\label{fig-3}
\end{figure}

\begin{figure}[ht]
  \centering
  \includegraphics[width=3in]{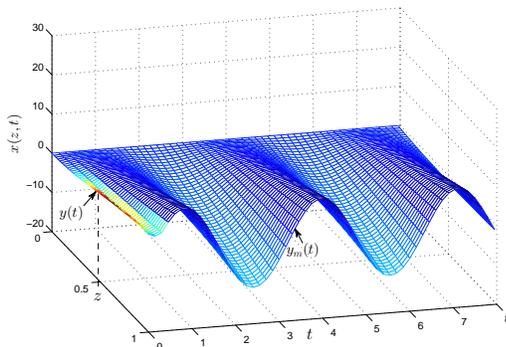}\vspace{-7mm}\\
  \caption{The evolution of the state $x(z,t)$ under the control of output feedback regulator (\ref{ym-contr-1})-(\ref{ym-contr-2})}\label{fig-4}
\end{figure}

It is easy to calculate ${C_{\Lambda m}}{\Pi _0}(z) = \left[ {\begin{array}{*{20}{c}}
{{\Pi _{01}}(1)}&{{\Pi _{02}}(1)}
\end{array}} \right] = \left[ {\begin{array}{*{20}{c}}
{ - 0.1229}&{ - 0.0868}
\end{array}} \right]$. Then, one can easily find $L_y$, e.g. ${L_y} = \left[ {\begin{array}{*{20}{c}}
0.1\\
{1}
\end{array}} \right]$ such that $S+L_yC_{\Lambda m}\Pi_0$ is exponentially stable. Consequently, the parameter of regulator $k$ can be computed through the explicit expression:
\[\begin{array}{l}
k(z) = {k_2}(z) = \frac{{(\Pi (z) + {\Pi _0}(z)){L_y}}}{{b(z)}}\\
 = \frac{{0.1\left( {{\Pi _1}(z) + {\Pi _{01}}(z)} \right) + \left( {{\Pi _2}(z) + {\Pi _{02}}(z)} \right)}}{{{e^{-0.5z}}}}
\end{array}\]
Finally, with the initial condition: ${r_w}(0) = \left[ {\begin{array}{*{20}{c}}
{0.1}\\
{4.6}
\end{array}} \right]$, the output feedback regulator (\ref{ym-contr-1})-(\ref{ym-contr-2}) is established. The results are shown in Figure \ref{fig-3} and Figure \ref{fig-4}.
In Figure \ref{fig-3}, we can see that in spite of the existence of disturbance $d(t)=5cos(2t)$, the controlled output $y(t)$ tracks the reference signal $y_r(t)$ well under the control of the proposed output feedback regulator (\ref{ym-contr-1})-(\ref{ym-contr-2}). In Figure \ref{fig-4}, the evolution of the state $x(z,t)$ are plotted. Moreover, the locations of $y(t)$ and $y_m(t)$ are pointed out. Therefore, we conclude that in this part, an output feedback regulator is constructed with the measured output $y_m(t)$ as its input, such that the controlled output $y(t)$ of the plant tracks the reference signal $y_r(t)$ despite of the disturbance $d(t)$.\vspace{-3mm}

\section{Conclusion}
In this paper, the output regulation problem of a mono-tubular heat exchanger process is considered. Based on the fact that the considered system is exponentially stable when the input variable is identically zero, a type of finite-dimensional output feedback regulator is constructed and the parameters are easily configured by applying the separation principle. Notably, for the output feedback regulator, only the measurement $y_m(t)$ of the plant (not exosystem) is available. Moreover, the measurement $y_m(t)$ and the controlled output $y(t)$ are allowed to be different. Finally, the mono-tubular heat exchanger process with the specified parameters is utilized to verify the main results and the performance of the proposed regulator is shown via the computer simulation.\vspace{-3mm}

\bibliography{myreference}             % bib file to produce the bibliography

\end{document}